\documentclass[a4paper,12pt]{article}
\usepackage[utf8x]{inputenc}
\usepackage{graphicx}
\usepackage[english]{babel}
\usepackage{amsmath}
\usepackage{amsfonts}
\usepackage{amssymb}
\usepackage{amsthm}
\usepackage{amstext}
\usepackage{pstricks}
\usepackage{color}
\usepackage{dsfont}
\usepackage{cite}

\newtheorem{theorem}{Theorem}[]

\theoremstyle{definition}

\theoremstyle{remark}

\def \cA {\mathcal{A}}

\def \cC {\mathcal{C}}

\def \cQ {\mathcal{Q}}

\def \cZ {\mathcal{Z}}

\def \k {\kappa}
\def \l {\lambda}
\def \s {\sigma}

\def \N {\mathbb{N}}

\def \R {\mathbb{R}}

\def \lra {\longrightarrow}

\def \al {\mathcal{A}^\ell}
\def \nal {\N^{\mathcal{A}^\ell}}
\def \ul {\lbrace\,1,\dots,\ell\,\rbrace}
\def \zl {\lbrace\,0,\dots,\ell\,\rbrace}
\def \nl {\N^{\ell+1}}

\def \On {(Z_n)_{n\geq 0}}

\def \on {Z_n}

\def \Zn {(Z_n)_{n\geq 0}}

\def \exa {e^{-a}}

\begin{document}

\begin{center}
\begin{LARGE}
Distribution of the quasispecies\\
for a Galton--Watson process\\
\vspace{5 pt}on the sharp peak landscape
\end{LARGE}

\begin{large}
Joseba Dalmau

\vspace{-12pt}
Universit\'e Paris Sud and ENS Paris

\vspace{4pt}
\today
\end{large}

\end{center}

\begin{abstract}
\noindent
We study a classical multitype Galton--Watson process
with mutation and selection.
The individuals are sequences of fixed length over a finite alphabet.
On the sharp peak fitness landscape 
together with independent mutations per locus,
we show that, as the length of the sequences goes to $\infty$
and the mutation probability goes to $0$,
the asymptotic relative frequency
of the sequences differing on $k$ digits from the master sequence
approaches 
$$(\s\exa-1)\frac{a^k}{k!}\sum_{i\geq1}\frac{i^k}{\s^i}\,,$$
where $\s$ is the selective advantage of the master sequence
and $a$ is the product of the length of the chains 
with the mutation probability.
The probability distribution $\cQ(\s,a)$ on the non negative integers
given by the above formula is the quasispecies distribution
with parameters $\s$ and $a$.
\end{abstract}

\section{Introduction}
Understanding the origin of life on Earth
has always been a major objective of science.
One of the many important contributions to the subject
is the 1971 article by Manfred Eigen~\cite{Eig},
which theoretically explores
the early stages of life on our planet.
As an explanation of how evolution of very simple macromolecules
might have occurred,
Eigen proposed a model known today as
Eigen's model or the quasispecies model.
The model aims at describing
the evolution of a population of macromolecules,
driven both by selection and mutation forces.
Due to the simplicity of the reproducing individuals,
Eigen's model can be synthesised 
as a system of differential equations,
obtained from the chemico--kinetic reactions
the different macromolecules are subject to:
$$x_k'(t)\,=\,
\sum_{i=1}^N f(i) Q(i,k) x_i(t)
-x_k(t)\sum_{i=1}^N f(i) x_i(t)\,,\qquad
1\,\leq\,k\,\leq\,N\,.$$
Here,
the different possible genotypes
are numbered from $1$ to $N$;
$x_k(t)$ represents the concentration
of individuals with genotype $k$ at time $t$;
$f(i)$ is the fitness (reproductive rate)
of the $i$--th genotype,
and $Q(i,k)$ is the probability
that an individual with genotype $i$
is transformed into an individual with genotype $k$
by mutation (unfaithful replication during reproduction).
Thus,
the first term in the differential equation
accounts for the production of genotype $k$ individuals,
while the second term 
accounts for the loss of individuals with genotype $k$;
the second term is proportional 
to the concentration of genotype $k$ individuals
as well as to the population's average fitness,
and it helps to keep the total concentration of chains constant.

One of the simplest scenarios we can consider
is that of the sharp peak landscape
together with independent mutations per locus.
In the sharp peak fitness landscape
all sequences but one,
the master sequence,
have the same fitness,
while the master sequence has a higher fitness than the rest.
Mutations happen during reproduction
independently on each locus of the sequence,
with equal probability.
Eigen studied this simple scenario
and found that two major phenomena take place.
The first is an error threshold phenomenon:
there is a critical mutation probability such that
for above--critical mutation probabilities
the equilibrium state of the population 
is a totally disordered one.
The second phenomenon 
is found for below--critical mutation probabilities:
in this case the equilibrium state of the population
is no longer disordered;
it contains a positive concentration of the master sequence,
together with a cloud of mutants
that closely resemble the master sequence.
This kind of distribution has come to be known as a quasispecies
distribution.

The concept of error threshold,
as well as that of quasispecies,
are very appealing to the scientific community,
mostly due to their potential 
for qualitatively explaining
the behaviour of a wide range of biological populations.
Since Eigen introduced them,
it has long been sought to extend the concepts
to many other situations,
both experimentally and theoretically.
From a theoretical point of view,
there are two main objections
to the applicability of Eigen's model
to more complex kinds of populations.
The first objection comes from considering at the same time
finite chain length and infinite population size:
if the individuals we seek to model are fairly complex,
the number of possible genotypes
largely exceeds the size of any viable population,
a feature that Eigen's model fails to account for.
The second objection is due to
the deterministic nature of Eigen's model:
again,
for fairly complex individuals,
the description of the reproduction mechanism
by chemico--kinetic reactions
is completely out of reach,
and it is typically replaced by some random mechanism.
The program is thus settled:
to retrieve the error threshold phenomenon
and a quasispecies distribution
for finite population stochastic models.
For a discussion 
on the several contributions to this program
we refer the reader to~\cite{CerfM,CD}.

In the series of papers~\cite{CerfM,CerfWF,CD,Dalmau},
the authors studied the classical 
Moran and Wright--Fisher models,
recovering both the error threshold phenomenon
and a quasispecies distribution
for mutation rates below the error threshold.
Furthermore, the quasispecies distribution
happens to be the same for both models,
and an explicit expression was found:
the concentration of sequences
differing in exactly $k$ digits 
from the master sequence is given by
$$(\s\exa-1)\frac{a^k}{k!}\sum_{i=1}^\infty \frac{i^k}{\s^i}\,,$$
where $\s>1$ is the reproductive advantage of the master sequence
and $a$ is the product of the mutation probability
with the length of the sequences.
We call this distribution
the quasispecies distribution with parameters $\s$ and $a$,
and we denote it by $\cQ(\s,a)$.
Both the Moran and the Wright--Fisher models
are constant population models,
since their aim is to describe 
a sufficiently large population 
which has stabilised in its environment.
However,
we might be interested in 
studying the evolution of a population
in its early stages.
The size of such a population 
is very likely to undergo significant fluctuations,
the classical stochastic model for this situation
is the Galton--Watson branching process.
The aim of our article
is to study a Galton--Watson branching process,
with selection and mutation,
in order to recover the phase transition phenomenon
and the quasispecies distribution.

Demetrius, Schuster and Sigmund~\cite{DSS}
already pursued this task
in a more general context:
a general fitness landscape
as well as a general mutation kernel.
In~\cite{ABCJ},
Antoneli, Bosco, Castro and Janini generalised the work in~\cite{DSS}
by studying a multivariate branching process,
which incorporates neutral, deleterious and beneficial mutations.
Our setting is closer to that of~\cite{DSS}
than~\cite{ABCJ};
our aim is to show that
for the sharp peak landscape
along with per--locus independent mutations,
the quasispecies distribution
is again the one obtained for the Moran model
and for the Wright--Fisher model.
In~\cite{DSS},
it was proved that 
the relative frequencies of the genotypes
converge to those given by 
the stationary solution 
of Eigen's system of differential equations.
However,
the quasispecies distribution
is a distribution on the Hamming classes of the sequence space,
which arises in a particular asymptotic regime.
Thus,
we cannot apply the results in~\cite{DSS} directly.
Along the lines of~\cite{CerfM,CerfWF,CD,Dalmau},
we develop our argument from scratch.
We start by
defining the Galton--Watson process on the genotypes,
with selection and mutation.
We formally show how to pass from the process on the genotypes
to a Galton--Watson process on the Hamming classes.
The relative frequencies of the classes
are shown to converge to the stationary solution
of the corresponding Eigen's system,
as shown in~\cite{DSS}.
Finally,
the stationary solution to this particular Eigen's system
is shown to converge to the quasispecies distribution.

Our article is organised as follows:
first we define a multitype Galton--Watson process
to model the evolution of a finite population.
We state next the main result of the article,
and all the remaining sections are devoted to the proof  
of the main result.

\section{The Galton--Watson process}
In this section
we define a multitype Galton--Watson process
driving the dynamics of a finite population,
which incorporates both selection and mutation effects.
Let us begin by introducing the individuals 
that will form the population.

\emph{Individuals.}
Let $\cA$ be a finite alphabet
of cardinality $\k\geq1$,
and consider sequences of fixed length $\ell\geq1$
over the alphabet $\cA$.
A sequence in $\cA^\ell$
represents the genotype of a haploid individual.
We study the evolution of a population
of such individuals,
with selection and mutation.

\emph{Sharp peak landscape.} 
The selection mechanism 
is given by a fitness function 
$A:\al\lra\R_+$.
Many fitness landscapes might be considered,
but we choose to work with the sharp peak landscape:
there is a particular sequence ${w^*\in\al}$,
called the master sequence,
whose fitness is $\s\geq1$,
while every other sequence in $\al$ 
has fitness $1$.
So, the fitness function in this case is given by
$$\forall u\in\al\qquad
A(u)\,=\,\begin{cases}
\quad\s&\quad\text{if }\ u=w^*\,,\\
\quad 1&\quad\text{if }\ u\neq w^*\,.
\end{cases}$$ 
\emph{Independent mutations per locus.} 
Mutations happen randomly
due to unfaithful replication of the chains,
independently on each locus of the chain,
with equal probability $q\in\,]0,1[\,$
for all loci.
When an allele mutates,
it does so to a randomly chosen letter, uniformly
from the $\k-1$ letters still available in the alphabet $\cA$.
This mutation mechanism
can be encoded into a mutation kernel in the following manner:
$$\forall u,v\in\al\qquad
M(u,v)\,=\,\prod_{i=1}^\ell
\Big(
(1-q)1_{u(i)=v(i)}
+\frac{q}{\k-1}1_{u(i)\neq v(i)}
\Big)\,.$$

\emph{The multitype Galton--Watson process}
is a Markov chain with values in 
$\smash{\N^{\k^\ell}}$,
$$X_n\,=\,
\big(X_n(u), 
u\in\al\big)\,,\qquad
n\geq 0\,.$$
For each $u\in\al$ and $n\geq 0$,
$X_n(u)$ represents
the number of individuals with genotype $u$
present in the population at time $n$.
At each generation,
each individual in the population
gives birth to a random number of children,
independently of the other individuals 
and of the past of the process.
The number of offspring of an individual $u\in\al$
is distributed as a Poisson random variable 
with mean $A(u)$.
The newborn individuals then mutate according to the kernel $M$.
The new generation is formed by all the offspring,
after mutation.

\emph{Generating functions.}
The classical tool for studying the Galton--Watson process
we just described is generating functions,
which are also useful for formally defining
the transition mechanism of the process.
Let $u\in\al$ and define the function
$f^u:[-1,1]^{\al}\lra\R$ by:
$$\forall s\in [-1,1]^{\al}\qquad
f^u(s)\,=\,
\sum_{r\in\N^{\al}}p^u(r)\prod_{v\in\al}s(v)^{r(v)}\,,$$
where $p^u(r)$
represents the probability
that an individual with genotype $u$
has $r(v)$ children with genotype $v$. For each $v\in\al$:
$$\forall r\in\N^{\al}\qquad
p^u(r)\,=\,
e^{-A(u)}A(u)^{|r|_1}
\prod_{v\in\al}\frac{M(u,v)^{r(v)}}{r(v)!}\,.$$
Here $|r|_1$ represents the usual 1--norm
of the vector $r$, that is,
the sum of its components.
For an initial population $X_0$
consisting of one genotype $u$ individual only,
$X_1$ is a random vector having generating function $f^u$.
In general,
for $n\geq0$,
if $X_n=r\in\nal$,
then $X_{n+1}$ is the sum of $|r|_1$ random vectors,
where, for each $u\in\al$,
$r(u)$ of the random vectors have generating function $f^u$.
Note that the null vector is an absorbing state.

\section{Main result}
Since we work with the sharp peak landscape fitness function,
we can classify the sequences in $\al$ 
according to the number of digits 
they differ from the master sequence.
Precisely,
the Hamming distance between two sequences $u,v\in\al$
is defined as the number of digits where the two sequences differ:
$$d_H(u,v)\,=\,\text{card}\big\lbrace\,
i\in\ul : u(i)\neq v(i)
\,\big\rbrace\,.$$
For each $k\in\zl$,
let $\cC_k$ be the set of the sequences
in $\al$ at Hamming distance $k$ from the master sequence:
$$\cC_k\,=\,\lbrace\,
u\in\al : d_H(u,w^*)=k
\,\rbrace\,.$$
We refer to the set $\cC_k$ as the $k$--th Hamming class.
Our aim is to study the concentration of the individuals of $X_n$
which are in the class $k$ in the following asymptotic regime:
$$\ell\to\infty\,,\qquad\quad
q\to 0\,,\qquad\quad
\ell q \to a\in [0,\infty]\,.$$
We have the following result.

\begin{theorem}
The process $(X_n)_{n\geq 0}$
has a positive probability of survival.
Conditioned on the event of non--extinction,
if $\s\exa\leq 1$ then
$$\forall k\geq 0\qquad
\lim_{\genfrac{}{}{0pt}{1}{\ell\to\infty,\,q\to0}{\ell q\to a}}
\,\lim_{n\to\infty}\,
\frac{1}{|X_n|_1}\sum_{u\in\cC_k}X_n(u)\,=\,0\,.$$
If $\s\exa>1$ then
$$\forall k\geq 0\qquad
\lim_{\genfrac{}{}{0pt}{1}{\ell\to\infty,\,q\to0}{\ell q\to a}}
\,\lim_{n\to\infty}\,
\frac{1}{|X_n|_1}\sum_{u\in\cC_k}X_n(u)\,=\,
(\s\exa-1)\frac{a^k}{k!}
\sum_{i\geq1}\frac{i^k}{\s^i}\,.$$
\end{theorem}
The right hand side in this equation
is the concentration of the $k$--th Hamming class
in the distribution of the quasispecies $\cQ(\s,a)$
with parameters $\s$ and $a$.
We devote the rest of the paper to the proof of this result.

\section{The occupancy process}
In this section we build an occupancy process
$$\On\,=\,(Z_n(0),\dots,Z_n(\ell))_{n\geq 0}\,,$$
to keep track of the number of sequences 
in each of the Hamming classes.
Here $Z_n(l)$ represents the number of individuals in $X_n$
that are at distance $l$ from the master sequence.
In order to build the occupancy process formally,
we use the classical lumping technique~\cite{KS}.

\emph{Fitness.} The fitness function $A$
can be factorised into Hamming classes:
define the function $A_H:\zl\lra\R_+$ by
$$\forall l\in\zl\qquad
A_H(l)\,=\,\begin{cases}
\quad\s&\quad\text{if }\ l=0\,,\\
\quad 1&\quad\text{if }\ 1\leq l\leq \ell\,.
\end{cases}$$
Then, for each $u\in\al$ we have
$A(u)=
A_H(d_H(u,w^*))$.

\emph{Mutations.} The mutation matrix $M$
can also be factorised into the Hamming classes.
Indeed, for each $u\in\al$ and $c\in\zl$, the value
$$\sum_{v\in\cC_c}M(u,v)$$
depends on $u$ through its Hamming class only
(lemma 6.1 in~\cite{CerfM}).
For $b,c\in\zl$, 
let us call $M_H(b,c)$ 
this common value for $u$ in $\cC_b$. 
The coefficient $M_H(b,c)$ 
can be analytically expressed as
$$
\sum_{
\genfrac{}{}{0pt}{1}{0\leq k\leq\ell-b}{
\genfrac{}{}{0pt}{1}
 {0\leq l\leq b}{k-l=c-b}
}
}
{ \binom{\ell-b}{k}}
{\binom{b}{l}}
q^k
(1-q)^{\ell-b-k}
\Big(\frac{q}{\kappa-1}\Big)^l
\Big(1-\frac{q}{\kappa-1}\Big)^{b-l}\,.
$$
\emph{Lumping.} 
Let $\cZ:\N^{\al}\lra\N^{\ell+1}$ 
be the map that associates to each population $r\in\N^{\al}$
the corresponding occupancy distribution:
$$\forall r\in\N^{\al}\quad
\forall l\in\zl\qquad
\cZ(r)(l)\,=\,
\sum_{u\in\cC_l}r(u)\,.$$
The occupancy process $\On$ is defined by
$$\forall n\geq0\qquad
\on\,=\,\cZ(X_n)\,.$$
We check next that the occupancy process 
is again a Galton--Watson process.
Let $k\in\zl$, $u\in\cC_k$, and $z\in\nl$.
We have
\begin{align*}
\sum_{\genfrac{}{}{0pt}{1}
{r\in\N^{\al}}{\cZ(r)=z}}p^u(r)\,&=\,
\sum_{\genfrac{}{}{0pt}{1}
{r\in\N^{\al}}{\cZ(r)=z}}e^{-A(u)}
A(u)^{|r|_1}
\prod_{v\in\al}\frac{M(u,v)^{r(v)}}{r(v)!}
\\&=\,
e^{-A_H(k)}A_H(k)^{|z|_1}
\sum_{\genfrac{}{}{0pt}{1}
{r\in\N^{\al}}{\cZ(r)=z}}\prod_{v\in\al}
\frac{M(u,v)^{r(v)}}{r(v)!}\,.
\end{align*}
Decomposing the last sum and product into Hamming classes,
we obtain
\begin{align*}
\sum_{\genfrac{}{}{0pt}{1}
{r\in\N^{\al}}{\cZ(r)=z}}\prod_{v\in\al}
\frac{M(u,v)^{r(v)}}{r(v)!}
\,&=\,
\prod_{l=0}^\ell\bigg(
\sum_{\genfrac{}{}{0pt}{1}
{t\in\N^{\cC_l}}{|t|_1=z(l)}}
\prod_{v\in\cC_l}\frac{M(u,v)^{t(v)}}{t(v)!}
\bigg)
\\&=\,
\prod_{l=0}^\ell \frac{M_H(k,l)^{z(l)}}{z(l)!}\,.
\end{align*}
Let 
$$p^k(z)\,=\,\sum_{\genfrac{}{}{0pt}{1}
{r\in\N^{\al}}{\cZ(r)=z}}p^u(r)\,=\,
e^{-A_H(k)}A_H(k)^{|z|_1}
\prod_{l=0}^\ell \frac{M_H(k,l)^{z(l)}}{z(l)!}\,.$$
Since this expression depends on $u$ only through $k$,
the sum
$$\sum_{\genfrac{}{}{0pt}{1}
{r'\in\nl}{\cZ(r')=z'}}
P\big(
X_{n+1}=z'\,|\,X_n=r
\big)$$
depends on $r$ only through $z=\cZ(r)$.
Thus,
by the classical lumping theorem,
the process $\On$ is a Markov chain
(the classical lumping theorem
is stated in~\cite{KS}
for finite state space Markov chains,
but both the result and the proof
carry over word by word to the case of 
denumerable Markov chains).
Let us define, for
$k\in\zl$ and $s\in [0,1]^{\ell+1}$,
$$f^k(s)\,=\,\sum_{z\in\nl}p^k(z)\prod_{l=0}^\ell s(l)^{z(l)}\,.$$
The process $\On$ is in fact a Galton--Watson process
with $\ell+1$ types,
having the following transition mechanism:
for all $n\geq 0$ and $z\in\nl$,
if $Z_n=z$,
then $Z_{n+1}$ is the sum of $|z|_1$
independent random vectors,
where,
for each $k\in\zl$,
$z(k)$ of the vectors have generating function $f^k$.

\section{Proof of theorem 1}
We use now the classical theory of branching processes~\cite{Harris}
in order to study the process $\Zn$.
The mean matrix $W$ of $\Zn$
is the matrix with coefficients 
$W(i,j),0\leq i,j\leq\ell$,
given by the expected number of class $j$
individuals in the first generation,
when the process starts with a population consisting
of just one individual in the class $i$.
The mean matrix is thus given by:
$$\forall i,j\in\zl\qquad
W(i,j)\,=\,A_H(i)M_H(i,j)\,.$$
The entries of the matrix $W$
are all positive.
By the Perron--Frobenius theorem,
there exist a unique largest eigenvalue $\l$ of $W$
and a unique positive and unitary eigenvector $\rho$
associated to $\l$. 
By the general theory of multitype Galton--Watson processes
(theorems~7.1 and~9.2 in chapter~II of~\cite{Harris}),
if $\l\leq 1$ then the population goes extinct with probability 1.
If $\l>1$ there is a positive probability of survival,
and conditioned on the event of non extinction, we have
$$\lim_{n\to\infty}\frac{Z_n(k)}{Z_n(0)+\dots+Z_n(\ell)}\,=\,
\rho(k)\,,\qquad
0\leq k\leq\ell\,.$$
From their definition, 
$\l$ and $\rho$ satisfy
$$\l\rho(k)\,=\,
\sum_{i=0}^\ell \rho(i)A_H(i)M_H(i,k)\,,\qquad
0\leq k\leq\ell\,.$$
Summing the above expression over $k$,
since $\rho$ is a unitary vector,
we deduce that
$$\l\,=\,
\sum_{i=0}^\ell\rho(i)A(i)\,=\,
(\s-1)\rho(0)+1\,.$$
Thus, the eigenvalue $\l$ is equal to the average fitness
of a population whose concentrations are given by the vector $\rho$.
We remark that solving the above system of equations
is equivalent to finding the stationary solutions
of the corresponding Eigen's system of differential equations.
From the above equation, we see that, in particular,
$\l\in\,]1,\s[\,$,
and this implies the first statement of the theorem:
the process $\Zn$ has a positive probability of survival.
It remains to study the asymptotic behaviour of $\l$ and $\rho$
when $\ell$ goes to $\infty$, 
$q$ goes to $0$ and $\ell q$ goes to $a$.
In this asymptotic regime 
the mutation kernel $M_H$ 
converges to the following limiting expression:
$$\forall i,k\geq0\qquad
\lim_{\genfrac{}{}{0pt}{1}{\ell\to\infty,\,q\to0}{\ell q\to a}}\,
M_H(i,k)\,=\,\begin{cases}
\quad \displaystyle\exa\frac{a^{k-i}}{(k-i)!} 
&\quad\text{if }\ k\geq i\,,\\
\quad 0
&\quad\text{if }\ k< i\,.
\end{cases}$$
Up to extraction of a subsequence, 
we can suppose that the following limits exist
$$\l^*\,=\,
\lim_{\genfrac{}{}{0pt}{1}{\ell\to\infty,\,q\to0}{\ell q\to a}}\,
\l\,,\qquad
\rho^*(k)\,=\,
\lim_{\genfrac{}{}{0pt}{1}{\ell\to\infty,\,q\to0}{\ell q\to a}}\,
\rho(k)\,,\quad k\geq 0\,.$$
Writing down the first equation of the system $\l\rho=\rho^T W$,
we see that
$$\s\rho M_H(0,0)\,<\,
\l\rho(0)\,<\,
\s\rho M_H(0,0)+
\max_{1\leq i\leq\ell}M_H(i,0)\,.$$
Since we also know that $\l>1$,
we conclude that $\l^*\geq\max\lbrace\,1,\s\exa\,\rbrace$.
As we have already pointed out,
$$\l\,=\,(\s-1)\rho(0)+1\,.$$
Thus, taking the limits in the above two equations
we deduce that
$$\l^*\,=\,(\s-1)\rho^*(0)+1\,,\qquad
\l^*\rho^*(0)\,=\,\s\rho^*(0)\exa\,.$$
Since $\l^*\geq\max\lbrace\,1,\s\exa\,\rbrace$,
we conclude that

$\bullet$ if $\s\exa\leq1$, then $\l^*=1$ and $\rho^*(0)=0$.

$\bullet$ if $\s\exa>1$, then 
$$\l^*\,=\,\s\exa\qquad
\text{and}\qquad
\rho^*(0)\,=\,\frac{\s\exa-1}{\s-1}\,.$$
Finally, writing down the $k$--th equation of the system
$\l\rho=\rho^T W$, we see that
\begin{multline*}
\s\rho_0M_H(0,k)+\sum_{i=1}^k\rho(i)M_H(i,k)
\,<\,\l\rho(k)\,<\,\\
\s\rho_0M_H(0,k)+\sum_{i=1}^k\rho(i)M_H(i,k)
+\max_{k<i\leq\ell}M_H(i,k)\,.
\end{multline*}
Thus, taking the limit we obtain the recurrence relation
$$\s\exa\rho^*(k)\,=\,
\s\rho^*(0)\exa\frac{a^k}{k!}
+\sum_{i=1}^k \rho^*(i)\exa\frac{a^{k-i}}{(k-i)!}\,,\qquad
k\geq1\,.$$
We conclude that if $\s\exa\leq1$, then $\rho^*(k)=0$ for all $k\geq0$,
and if $\s\exa>1$, then
$$\rho^*(k)\,=\,(\s\exa-1)\frac{a^k}{k!}
\sum_{i\geq1}\frac{i^k}{\s^i}\,,\qquad
k\geq0\,.$$
This can be seen by solving the recurrence relation by the method 
of generating functions, see, for example~\cite{CD}.

\bibliography{galwat}
\bibliographystyle{plain}
\end{document}